\documentclass[pre, %showpacs, 
floatfix, twocolumn]{revtex4}

\usepackage{graphicx}
\usepackage{amsmath}
\usepackage{latexsym,amsfonts} 
\usepackage{listings}
\usepackage{boxedminipage}
\usepackage{subfig}

\usepackage{floatflt}

\usepackage{array}

\usepackage{color}

\newcommand\cpp{$C$++}
\newcommand\cuda{CUDA}

\newcommand\cproglang{$C$}

\newcommand\openmp{OpenMP}

\newcommand\osp{OSP}
\newcommand\pmask{PM}
\newcommand{\prp}{PRP}
\newcommand\ps{PS}
\newcommand{\foutra}{Fourier transform}
\newcommand{\maproj}{method of alternating projections}

\def\mybeginwidetext{\begin{widetext}}
\def\myendwidetext{\end{widetext}}

\definecolor{lgrey}{gray}{.9}
\begin{document}

\title{ Real-Time Phase Masks for  Interactive 
 Stimulation of Optogenetic Neurons 
}

\author{Stephan C. Kramer,$^{1}$, Johannes Hagemann,$^{2}$ and D. Russell Luke$^{1}$}

\address{
$^1$ Institut f\"ur Numerische und Angewandte Mathematik, Georg-August-Universit\"at G\"ottingen
\\
$^2$ Institut f\"ur R\"ontgenphysik, Georg-August-Universit\"at G\"ottingen
}

\email{stkramer@math.uni-goettingen.de} 

\begin{abstract}
Experiments with networks of optogenetically altered neurons require stimulation 
with high spatio-temporal selectivity.
Computer-assisted holography is an energy-efficient method for 
robust and reliable 
addressing of single neurons on the millisecond-timescale inherent to biologial information processing.
We show that real-time control of neurons can be achieved by a CUDA-based hologram computation.
\end{abstract}

\maketitle
%\ocis{070.6120, 090.1760, 090.1995, 090.5694, 100.5070, 100.3190}
\section{Introduction}
\label{sec:Introduction-ph}

The development of light-sensitive neurons has been a
milestone in optogenetics~\cite{Deisseroth2006}.
The ability to engineer a neuron's optical sensitivity by genetic manipulation
is crucial for a non-destructive and
fast, yet accurate 
photostimulation (PS) of individual sites in networks of living neurons.
In vivo interaction with individual neurons is fundamental for a concise experimental study of such basic neurological processes like the mechanisms of learning. 
In terms of energy efficiency and spatial resolution
holographic methods are considered to be the most suitable for {\ps}~\cite{Golan2009a}. 

\begin{figure}[pb]
\centering
\includegraphics[ 
height=5.8cm
]{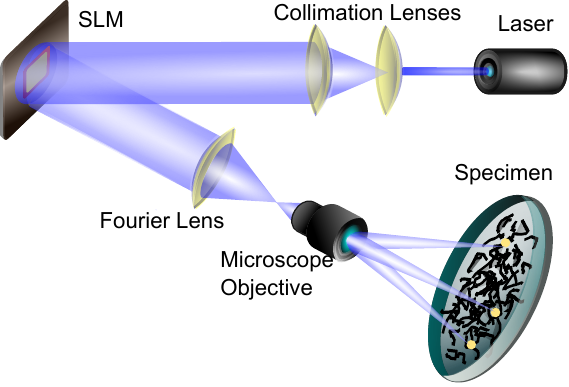}
\caption{Holographic illumination of a network of optogenetically altered neurons.
\vspace{-.2cm}
}
\label{fig:exp-setup}
\end{figure}

Holograms, i.e. 
computer generated phase masks (PM), 
displayed on a spatial light modulator (SLM)
realize pixel-wise phase retardations of a coherent laser beam.
Upon illumination the intensity of the Fourier transform of the {\pmask} yields a high-resolution optical
stimulation pattern 
(OSP) 
at the specimen. For a sketch of the experimental setup see Fig.~\ref{fig:exp-setup}.
The {\osp} follows from the subset of neurons selected for stimulation.
By targeting specific neurons the neural activity in the network and thus its collective behavior can be influenced.

The basis of neural activity is the generation of spikes in the membrane potential at the axon hillock due to synaptic input.
The spikes travel along the axon to 
the synaptic connections to other nerve cells.
In genetically altered neurons light-sensitive ion channels are expressed in the cell membrane.
If lit with the correct frequency the ion channels open and thus change the membrane potential.
This either inhibits or enhances spike generation.
After transmission to the next neuron the spike adds to the synaptic input which may lead to another spike.

For interactive modification of the spiking behavior the optical stimulus must be generated within the time a spike 
needs to travel from one neuron to another.
Interspike intervals of adjacent neurons are in the range of 10-20 ms
and set the time-scale for computing the unknown \pmask.
Due to this severe time constraint multi-site stimulation thus had to use precomputed {\pmask}s, up to now.
For interactive network control {\pmask}s must be computed online which for 
frame rates in the required range of 0.1 to 1~kHz poses a substantial challenge.
On current many- and multi-core processing units this requires extensive parallelization.

Mathematically, computing a {\pmask} for a
given {\osp} constitutes an inverse problem equivalent to 
wavefront reconstruction (see \cite{Luke2002} and references therein)
and
is an instance of the phase retrieval problem (PRP) 
in diffraction imaging~\cite{Rayleigh92}.
Numerical approaches to the phase problem 
abound, but convergence results 
and global solutions are limited to special cases \cite{Hauptman, BCL1} that 
do not 
necessarily
apply to the case discussed in this paper. 
An arbitrary {\osp} is unlikely to have a phase-only
{\foutra}.
Thus our {\prp} is fundamentally 
{\em inconsistent} as defined in~\cite{Luke11}.  
To account for the mathematical structure, a careful analysis of the performance of the parallelization techniques available and a strong focus on long-term software reusability distinguishes this work from others, e.g.~\cite{Masuda:06,Shimobaba:10,Weng:12}.
Useful approximations of a {\pmask} for a given {\osp} can be achieved by
iterative algorithms like the widely used {\em Method of Alternating Projections}~\cite{Neumann50}, also known as~Gerchberg-Saxton algorithm \cite{Gerchberg1972}.

In this work 
we will combine parallel computation on graphics cards
with C++-based generic programming and 
a sound mathematical theory.
Only this combination of techniques allows to generate phase masks within~less than 10~ms, matching the dynamics of neural activity.

\section{Method of Alternating Projections} \label{sec:algs-impl}

\noindent
The wavefront is to be altered by a phase shift at the finite number of pixels of the SLM. The entire
system is modeled on a finite dimensional vector space. 
Let~$L_x, L_y>0$ be the dimensions of the SLM and $n_x, n_y$ the respective number of pixels.
We seek a signal~$u\in \mathbb{C}^N$ for~$N\equiv n_x\times n_y$.
The intensity distribution of the laser beam sets the amplitude of the wavefront~$u$ on the SLM.
Assuming a
constant intensity over one pixel, we discretize the intensity distribution by the nonnegative~$p\in \mathbb{R}^N$.  
Wavefronts~$u$ emanating from the SLM are 
given by the set of vectors \begin{eqnarray}
\label{eq:slm-constraint}
S & \equiv & \{ u \in \mathbb{C}^N ~: ~|u_{jk}| = p_{jk}, \nonumber \\
&& \quad j=1,2,\dots,n_x, ~k=1,2,\dots,n_y \}. 
\end{eqnarray}
Propagation of the light through the lens system is modeled by Fraunhofer diffraction~\cite{Goodman}.
The light at the SLM is related to the observed {\osp} at the
specimen by  the Fourier transform $F$. 
Waves with modulus matching
the amplitude distribution~$m\in \mathbb{R}^N$ of the {\osp} form
the set 
\begin{eqnarray}
\label{eq:pattern-constraint}
M \equiv \{ u \in  \mathbb{C}^{N}  ~: ~|(F u)_{jk}| = m_{jk}, \quad j \le n_x,  ~k \le n_y
\}.
\end{eqnarray}
On the one hand our wavefront must fulfill the amplitude constraint~Eq.~(\ref{eq:slm-constraint}),
on the other hand the amplitudes of its Fourier transform are fixed by the intensity distribution of the 
stimulation pattern,~Eq.~(\ref{eq:pattern-constraint}).
Altogether, the mathematical problem we address is to 
\begin{equation}\label{e:feas}
\mbox{Find}~u\in S\cap M.
\end{equation}
For a nonempty intersection the problem is defined to be {\em consistent};
otherwise {\em inconsistent} or {\em ill-posed}.  
A common algorithm for problems of this type is
the {\maproj}~\cite{Neumann50, Gerchberg1972}.
For a review of this and other projection-based approaches for the {\prp} 
see \cite{BCL1}.  Algorithms of this kind are built on {\em projection operators}
onto the constraint sets $S$ and $M$.  
By a {\em projection} of a point $u$ in a space $X$ onto a subset $C$ of that space, we mean 
the mapping of that point to the set of nearest points in $C$ with respect to 
the norm induced by the real inner product on $X$.
For general {\prp}s, it was proved in~\cite{Luke2002} 
that \begin{subequations}
\label{eq:projections}
\begin{eqnarray}
\label{e:PS}
\hspace{-.05cm}
P_Su  & \hspace{-.15cm}= \hspace{-.15cm}& \left\{v\,:\, v_{jk}=\begin{cases}
p_{jk} \frac{u_{jk}}{|u_{jk}|}, & \text{if
$u_{jk} \neq 0$;} \\
p_{jk}\exp(i\theta), & \mbox{for } \theta\in[0,2\pi)\end{cases}\right\},  \hspace{.7cm}  \\
\label{e:PM}
\hspace{-.05cm}P_Mu &\hspace{-.15cm}= \hspace{-.15cm}& \left\{F^{-1}\widehat u~:~\widehat u\in \widehat M(u)\right\}   
\end{eqnarray}
are projections onto the sets 
$S$ and~$M$, respectively, where
\begin{equation}
\label{e:Mhat}
\widehat M(u) \equiv \left\{\hat{v}~:~\hat v_{jk}
= 
\begin{cases}
m_{jk} \frac{(Fu)_{jk}}{|(F u)_{jk}|}, & \text{if }(Fu)_{jk} \neq 0 \\
m_{jk}\exp(i\theta), & \mbox{for } \theta\in[0,2\pi)\end{cases}\right\}.
\end{equation}
\end{subequations}
For given $u^0\in \mathbb{C}^N$ the {\maproj}
computes the iterates $u^{\nu}$ via
\begin{equation}\label{e:AP}
u^{\nu+1}\in P_SP_Mu^\nu, \quad \nu=0,1,2,\dots
\end{equation}
The multi-valuedness of Eq.~(\ref{eq:projections}) makes Eq.~(\ref{e:feas}) a {\em non-convex} feasibility problem 
\cite{Luke2002}.
Hence Eq.~(\ref{e:AP}) must be understood as a selection from set-valued mappings.  
Due to nonconvexity, except in special cases \cite{Hauptman}, global convergence
of Eq.~(\ref{e:AP}) cannot be guaranteed in general.  
For consistent {\prp}s local convergence results are available \cite{Luke11}.
Yet, it is more the exception than the rule that our {\prp} will be consistent:  
a set of fixed amplitudes cannot produce an arbitrary~{\osp}.  
Our 
numerical 
experiments
indicate 
that our {\prp}s are indeed inconsistent as 
measured by the magnitude of the {\em gap} 
\begin{eqnarray}
\label{eq:gap-error}
G & \equiv & \| P_{S}u^{\nu} - P_{M}u^{\nu}   \|_2 
\end{eqnarray}
between accumulation points in~$M$ and 
their projections onto~$S$.
The gap is measured in the standard Euclidean norm~$\|\cdot \|_2$.
This systematic inconsistency is a major difference between 
optogenetic {\ps} and {\prp}s due to imaging experiments. In the latter the diffraction pattern 
comprising the set $M$ is {\em causal}, 
that is, comes from diffraction by a physical object, e.g. a protein crystal.
Assuming that Eq.~(\ref{e:feas}) is inconsistent, we content ourselves with 
finding {\em best approximation pairs} $(u^*,v^*)\in \mathbb{C}^N\times\mathbb{C}^N$
with $u^*\in S$, $v^*\in M$, $P_Mu^*=v^*$ and $P_Sv^*=u^*$.  

To account for the particularities of
optogenetic {\ps}, we define the physical error as the sum
of the pixel-wise relative violation of deviation tolerances between target and
reconstructed {\osp}~\cite{HeckesError}. We allow for a relative deviation $t_{\ell} =0.1$ for non-zero target
pixels and an absolute deviation of  $t_{\text{d}}=3\cdot 10^{-4}$ for non-lit
pixels.
Violations are summed
in multiples of $t_{\ell}$ and $t_d$. With $u_{jk}^{\nu}$ being the intensity from the current iteration step, $m_{jk}$ the intensity in the target {\osp} and $\Theta(\cdot)$ the Heaviside step function, the total error is the sum of
\begin{subequations}
\label{eq:phys-error}
\begin{eqnarray}
\label{eq:phys-error-a}
E_\ell^{\nu}&= & \hspace{-.25cm}  \sum\limits_{pixels \in lit} \hspace{-.1cm}   
\left( \frac{ t_{\text{d}} |  m_{jk} - u_{jk}^{\nu} |}{t_\ell  m_{jk}}-t_{\text{d}}  \right) \nonumber \\
&& \times
\Theta\left( \frac{|  m_{jk} - u_{jk}^{\nu} |}{ m_{jk} }-t_\ell \right), 
\quad ~~ \\
\label{eq:phys-error-b}
E_{\text{d}}^{\nu}&= &  \hspace{-.25cm}   \sum\limits_{pixels \in dark} \hspace{-.1cm}   \left( u_{jk}^{\nu} - t_{\text{d}} \right) \cdot\Theta\left( u_{jk}^{\nu} - t_{\text{d}} \right)\,.
\end{eqnarray}
\end{subequations}

\section{Unified Implementation}
\label{phase-holos-unified-impl}

\mybeginwidetext
~
\begin{figure}
\subfloat[]{
\includegraphics[width=.475\textwidth, clip=true, trim=32mm 31.5mm 94mm
20mm]{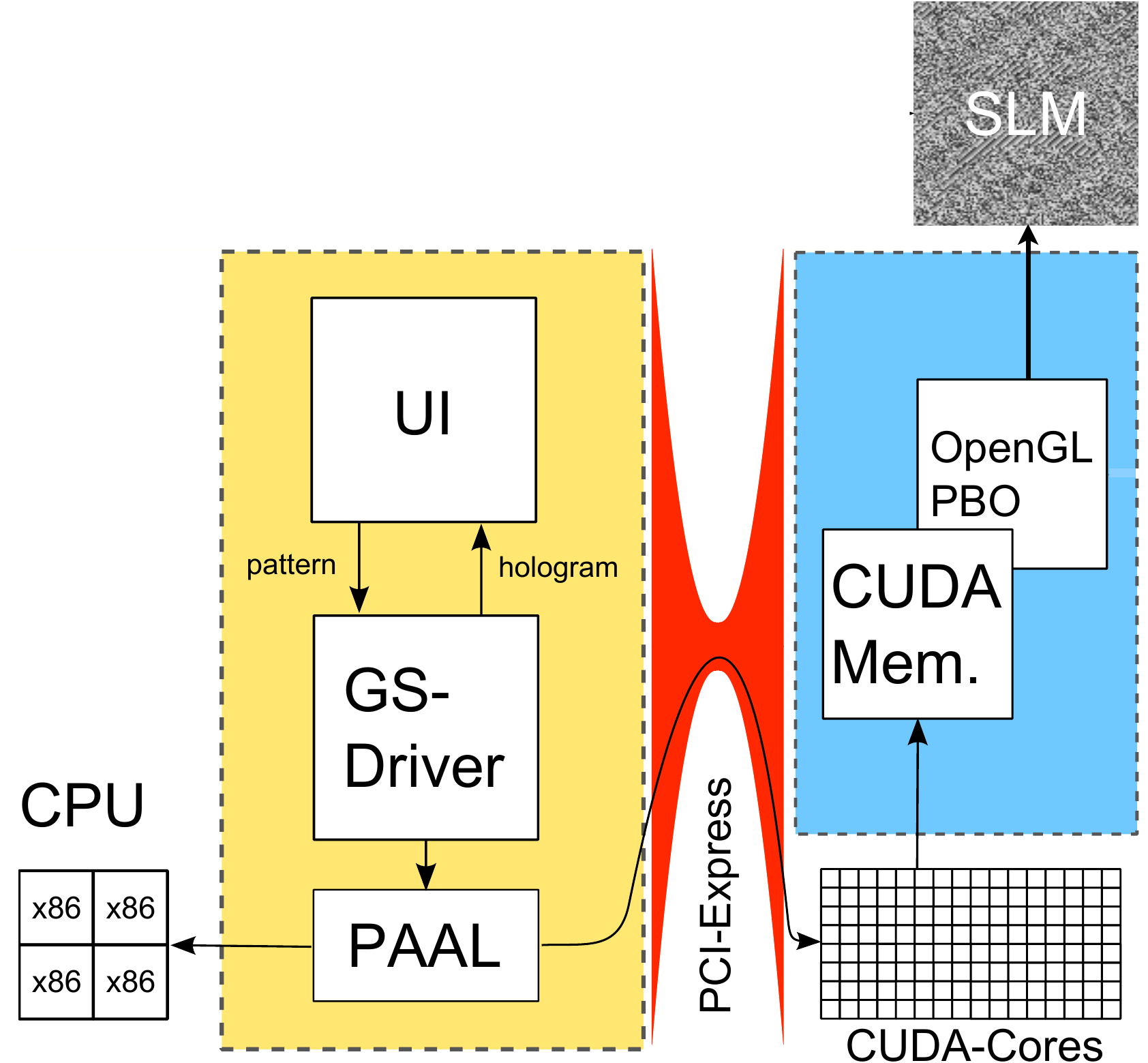}
\label{fig:program}
}
\subfloat[]{
\includegraphics[width=.48\textwidth, clip=true, trim=90mm 29mm 68mm 52mm
]{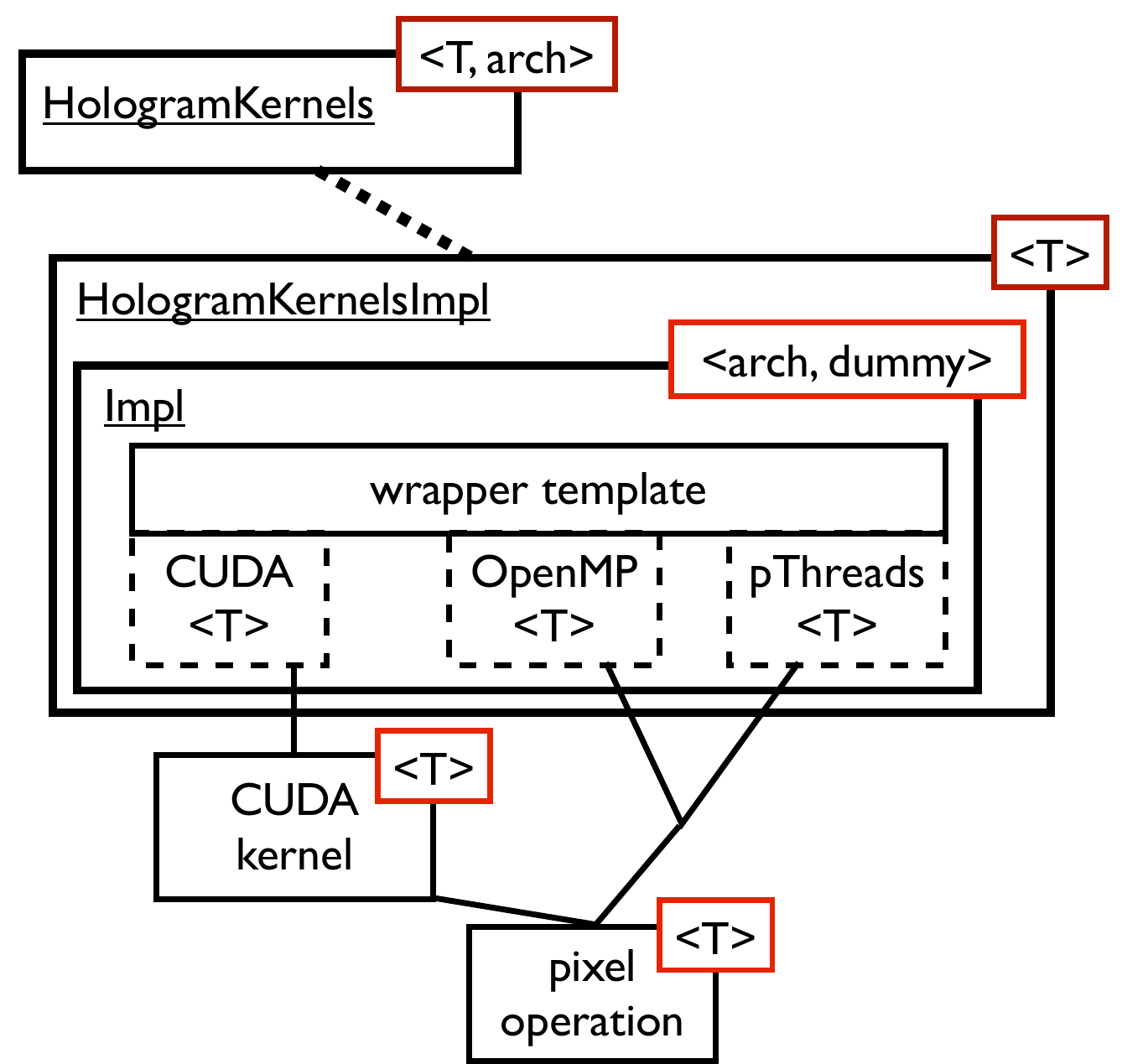}
\label{fig:PAAL-class-structure}
}
\caption{
Left: Software structure and its association with the different hardware components.
Right: Class diagram of the PAAL concept. Class names are underlined. Template arguments are given in red boxes. Dashed boxes indicate partial template specializations.}
\end{figure}
\myendwidetext
\noindent 
The compound system of CPU and GPU, each with dedicated memories, represents a non-uniform memory access architecture with a very heterogeneous distribution of processing capabilities and internal transfer rates.
Figure~\ref{fig:program} sketches
the class structure
and its distribution over the compound system of CPU and GPU.
The major bottleneck is the PCIe-BUS.
According to the PCIe v2.0 specification it has a maximal transfer rate of 8 GByte/sec although in practice one rather gets 4~to~5 GByte/sec.
This will rise to 16 GByte/sec with the forthcoming PCIe v3.0, yet memory transfer rates within the GPU are of the order of 100~GByte/sec.
On CPUs with integrated memory controllers the transfer rates are in the range of 20~to~30~GByte/sec.
To anticipate the rapid evolution of hardware and parallelization techniques we spent considerable effort on modularizing the program using C++'s templating capabilities. 
CUDA extends~{\cproglang} for programming NVidia GPUs.
OpenMP provides multithreaded parallelization on multi-core~CPUs.
Depending on the parallelization technique the program works on different sides of the PCIe-BUS.
To separate hardware-specific optimizations at the low-level, e.g. the details of Eq.~(\ref{eq:projections}), from the implementation of algorithms we introduced the concept of a parallel architecture abstraction layer~(PAAL).
Since most of our computational tasks are data-parallel they are perfect candidates for abstraction with respect 
to floating-point precision and parallelization. 
This is achieved by a suitable set  of 
template parameters leaving the generation of the hardware-specific part of the code to the compiler.
As FFTs in the projector onto the constraint representing the {\osp}, Eq.~(\ref{e:PM}), we use either  NVidia's cuFFT or the FFTW~\cite{Frigo2005} which offers an OpenMP- as well as a pthreads-based variant.

The aim of the PAAL concept is a quick recombination of algorithms and
parallelization strategies by explicit template specializations.
The front end comprises the user-interface (UI) and manages the execution of the phase retrieval algorithms for which separate driver classes exist, e.g. GS-DRIVER for the {\maproj}.
The final {\pmask} is transferred to an OpenGL framebuffer object
for display on the SLM, cf. Fig.~\ref{fig:program}.

The PAAL concept is explained best by walking through the essential parts of its class structure.
This is done roughly in a top down approach, i.e. from host to device and how things build on each other.
A sketch of the class structure is given in Fig.~\ref{fig:PAAL-class-structure}.
At the top is the interface to the GS-Driver which is formed by the class \lstinline{HologramKernels}. It takes two template arguments:  \lstinline{T} for the precision and  \lstinline{arch} for the architecture the algorithm is to run on.
At the bottom of the hierarchy is the operation one has to do on a particular pixel. 

To express that the class \lstinline{HologramKernels} is {\it implemented with} \lstinline{HologramKernelsImpl} inheritance is private~\cite{meyers2005effective} (indicted by the dashed line in Fig.~\ref{fig:PAAL-class-structure}). 
The class \lstinline{HologramKernels} needs partial specializations for the different architectures because 
for the {\cuda} kernels the wrapper functions behave differently with respect to the architecture. 
On a NVidia GPU they have to call a {\cuda} kernel. 
On a CPU they have to either use {\openmp} or pthreads for parallelization.
The parallel execution of the per-pixel operation via {\openmp} or pthreads can be done directly in the specialization of the wrapper function. In the following we omit the pthreads specialization as it is structurally very similar to {\openmp}.
The wrapper functions are implemented by the internal class \lstinline{Impl} of \lstinline{HologramKernelsImpl}.
The reason for this particular design is that the {\cpp} standard does not allow to define partial specializations 
of (a subset of) the member functions of a class.
This issue can be circumvented by introducing an internal class with a dummy template parameter and to partially specialize its members.
Within the class \lstinline{Impl} the particular type of real and complex numbers is deduced from the template parameter~\lstinline{T} by means of 
a suitable \lstinline{Traits} 
structure. 
This is a typical example of template metaprogramming~\cite{AbrahamsGurtovoy2004}.

The back end, i.e. the details of the implementation, are stored in a separate source file to keep g++ away from \cuda-specific code.
Within the evaluation of Eqs.~(\ref{e:PS})~and~(\ref{e:Mhat}) we need precision-dependent tolerances for what is considered as zero. 
To this end we localize the inevitable magic numbers in a structure~\lstinline{__eps} and a function~\lstinline{__is_zero}.
In case of being compiled into a {\cuda} kernel the \lstinline{__device__} keyword is put into effect indicating that the function can only be executed on the device, i.e. the GPU.
Prepending  \lstinline{__device__} by \lstinline{__host__} signals the compiler (i.e. nvcc) to compile two versions of a function or operator. 
One for the execution on the GPU and one for the CPU.
At the binary level these are distinct functions.

The actual per-pixel operation is done by an architecture-independent function~\lstinline{__ps_element}.
Its arguments are a pointer \lstinline{d_devPtr*} to the beginning of the array of pixels of the iterated image~$Fu^{\nu}$, a pointer \lstinline{d_original*} to the beginning of the array of pixels of the original image and the lexicographic index of the pixel \lstinline{x}.
The  {\cuda}  kernel~\lstinline{__ps}~(listing~\ref{lst:ph-cuda-kernel}) basically has the same arguments as the element function. 
The kernel gets the size of the image as additional argument in order to avoid operating on non-existent pixels.
The kernel computes the position~\lstinline{x} of its pixel from its \lstinline{threadIdx} and \lstinline{BlockIdx}. Given the pixel position is within the bounds of the array the element function is invoked.
\begin{lstlisting}[caption={PAAL: CUDA kernel for amplitude adaption}, label={lst:ph-cuda-kernel}]
template <typename T>
__global__ void __ps(T *d_devPtr, const T *d_original, const int size)
{
int x = blockDim.x*blockIdx.x + threadIdx.x;

if(x < size)
__ps_element<T, gpu_cuda>(d_devPtr, d_original, x);
}
\end{lstlisting}

The missing link between back end and driver class is the specializations of the wrapper functions for the kernels.
The GPU version (listing~\ref{lst:ph-cuda-wrapper}) starts as many threads as there are pixels for the kernel~\verb+__ps+. 
Each thread starts the device function \verb+__ps_element+, so that each pixel (vector element) is processed.
\begin{lstlisting}[caption={PAAL: GPU specialization of wrapper function}, label={lst:ph-cuda-wrapper}]
template<typename T>
template< typename dummy>
void 
HologramKernelsImpl<T>::Impl<gpu_cuda, dummy>::ps
(Complex *d_devPtr, const Complex *d_original, const int size)
{
int threads_per_block = 512;
int blocks = (size + threads_per_block - 1) / threads_per_block;

__ps<T><<<blocks, threads_per_block>>>(d_devPtr, d_original, size);

cudaThreadSynchronize();
}
\end{lstlisting}
The CPU-{\openmp} version (listing~\ref{lst:ph-openmp-wrapper})  has a ~\lstinline{for}-loop over all pixels in the image
which is parallelized by an {\openmp} preprocessor directive.
By declaring the \verb+__ps_element+ to be a \verb+__host__ __device__+
function, we can call the same function from the CPU as from the GPU but without the intermediate kernel layer.
In this way we have
the actual computation implemented only once. With the individual specialized classes wrapped
around this implementation we can choose our computing precision and hardware.

\begin{lstlisting}[caption={PAAL: CPU specialization of wrapper function}, label={lst:wrapper-openmp}, label={lst:ph-openmp-wrapper}]
template<typename T>
template< typename dummy>
void 
HologramKernelsImpl<T>::Impl<cpu, dummy>::ps
(Complex *d_devPtr, const Complex *d_original, const int size)
{
#pragma omp parallel for private(i)
	for(int i = 0; i < size; i++)
		__ps_element<T, cpu>(d_devPtr, d_original, i);
}
\end{lstlisting}
Finally, we have to provide full template specializations of all the combinations of precision and architecture template 
parameters we want to work with.
This must be at the end of the hardware-specific source file as all functions have to be declared and their bodies defined before the
class can be explicitly instantiated by the compiler~\cite{VandervordeJosuttis2002}.
The explicit specializations are necessary since we compile the back end with a different compiler than the front end of the program.
\mybeginwidetext
~
\begin{figure}[bp]
\centering
\includegraphics[angle=0,  width=.77\textwidth, clip=true, trim = 32mm 103mm 174.5mm 20mm]
{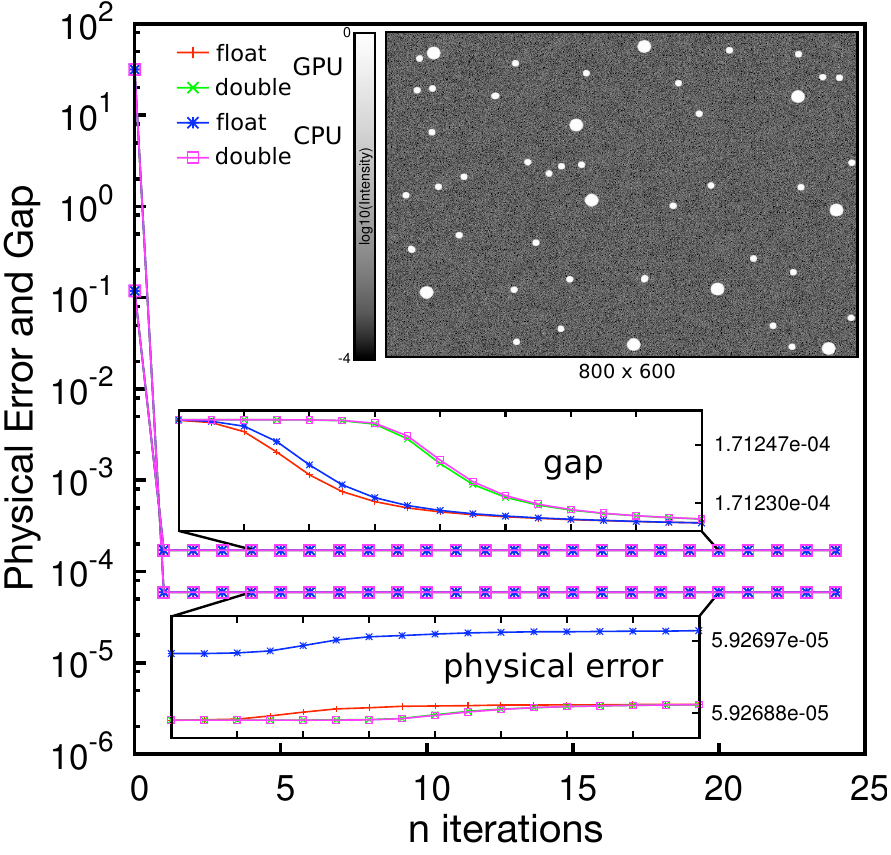}

\caption{Convergence for a spot pattern of the physical error and the constraint gap.
}
\label{fig:results}
\end{figure}
\myendwidetext

\mybeginwidetext
~
\begin{figure}[htbp]
\centering
\includegraphics[angle=0, width=.77\textwidth, 
clip=true, trim = 165mm 103mm 42mm 20mm]
{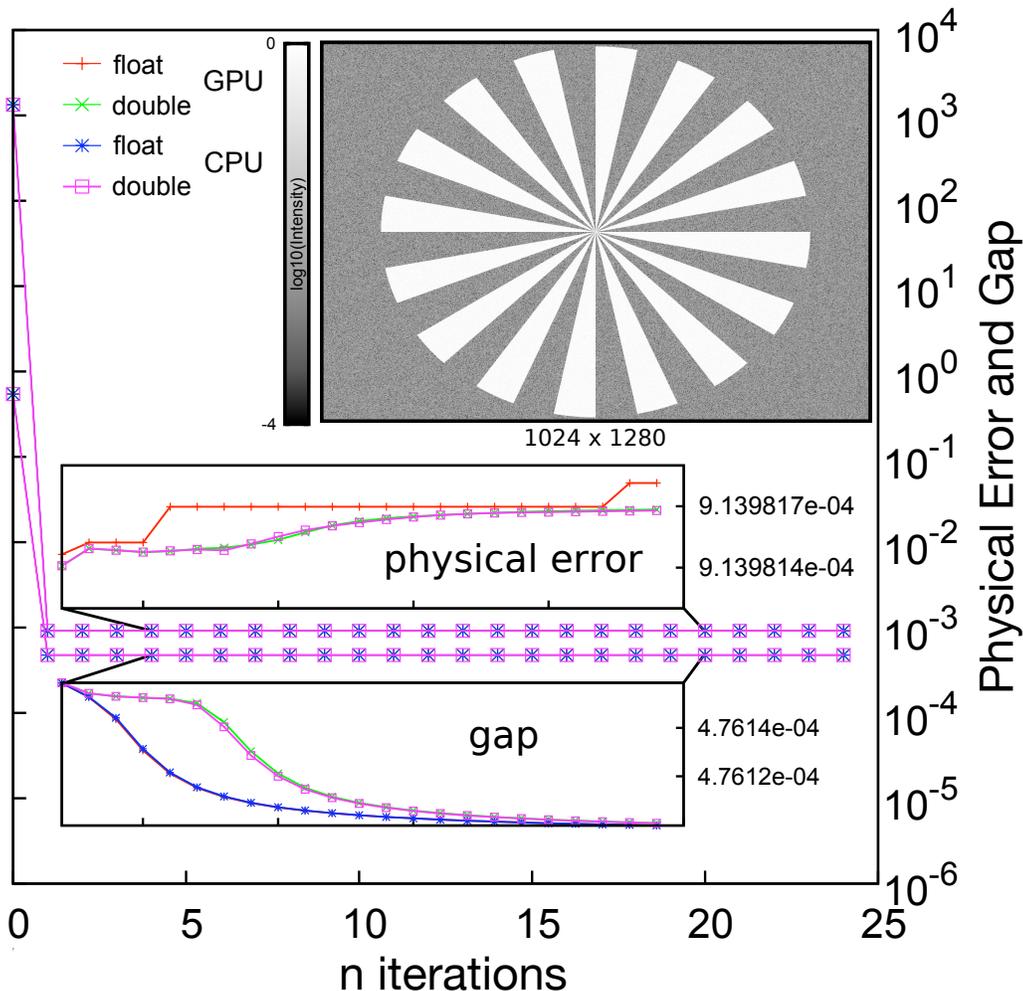}
\vspace{-.15cm}
\caption{Siemens star: Convergence 
for single and double precision.}
\label{fig:results-siemens}
\end{figure}
\myendwidetext
\section{Results}
\label{sec:results-ph}
The {\osp} for benchmarking the computation of the phase masks for photo-stimulation are bright spots on a dark background (Fig.~\ref{fig:results}).
The limits of spatial resolution in the reconstruction is tested on the Siemens star~(Fig.~\ref{fig:results-siemens}).  
In both cases we use $P_{M}\hat u^{0}$ as initial {\pmask} where $u^{0}$ is the 1-bit target {\osp}.
Thus, our initial condition is computed by taking the Fourier transform of the {\osp}, adapting the Fourier coefficients to the amplitude constraints on the SLM and transforming back into real space again.
The other obvious choice as initial condition would be to use random phases. 
The resulting phase masks do not differ significantly from those obtained by using $P_{M}\hat u^{0}$ but take longer to converge.
Therefore we skip them in the following discussion.
\subsection{Benchmarking}
For interactive holographic {\ps}  the
physical error,~Eq.~(\ref{eq:phys-error}), must reach a sub-threshold level within interspike intervals, i.e. 10 to 20~ms.
Hence, 
the first issue is which parallelization technique provides sufficient performance to meet this requirement.
The GPU-based computations were done on a Tesla C2070. 
The CPU-based ones using OpenMP or pthreads were run on a two-socket system with X5650 Xeons.
The CPUs have six cores each. Therefore we decided to use 12 threads, i.e. as much as there are physical cores.

\begin{figure}[hbp]
\centering
\subfloat
{
\includegraphics[width=.401\textwidth , clip=true,  trim = 122.5mm 106.2mm 132mm 62.25mm 
]
{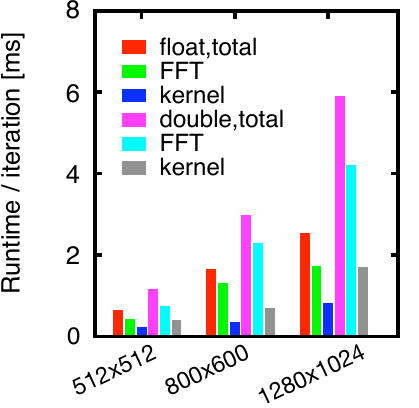}
\label{fig:runtimes}
}
\subfloat
{
\includegraphics[width=.4\textwidth, clip=true,  trim = 24mm 73.5mm 133mm 21.0mm]{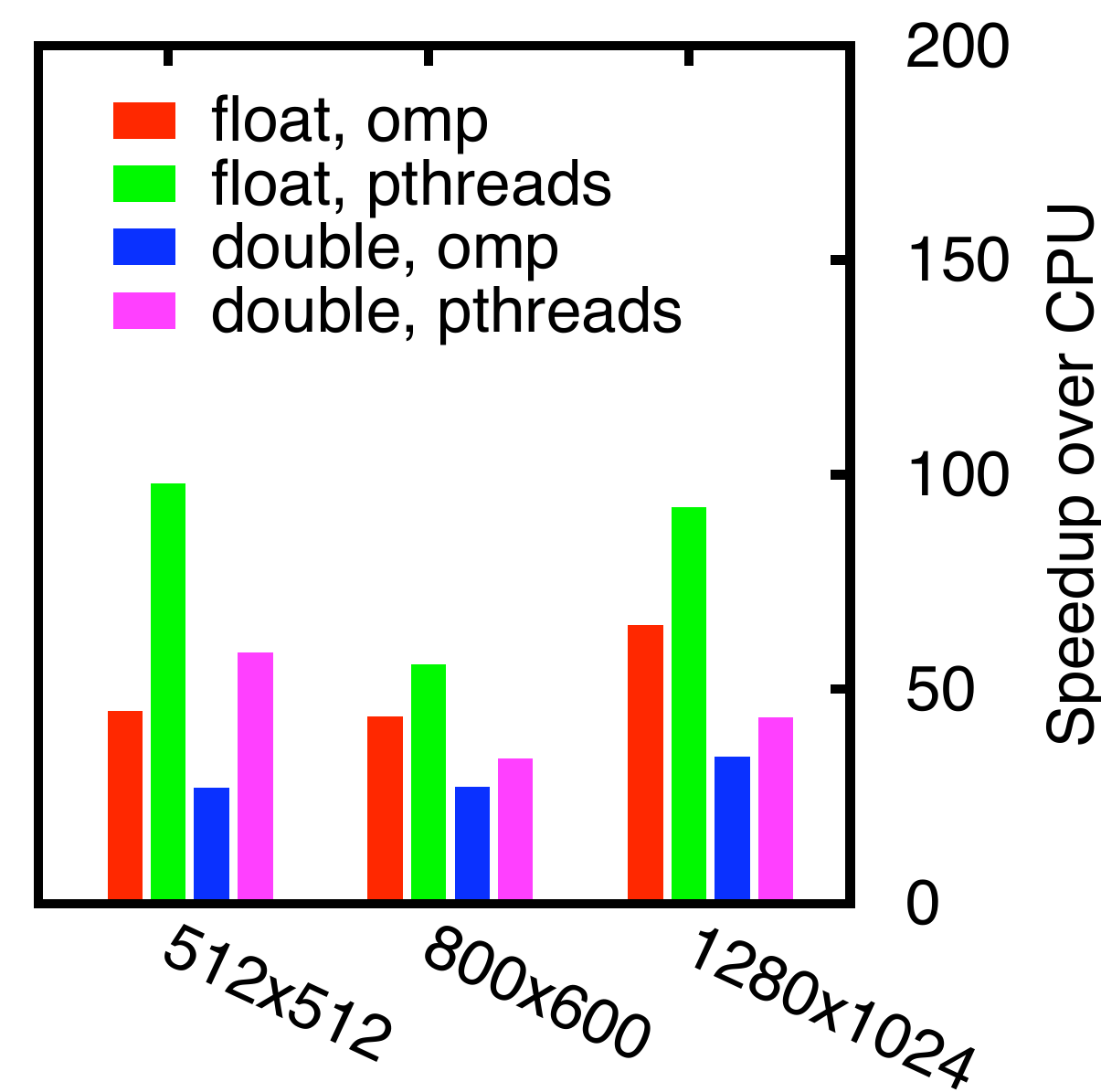}
\label{fig:speedup} 
}
\caption{Left: GPU runtimes per iteration for various image sizes and 25 iterations of Eq.~(\ref{e:AP}). Right: Speedup per iteration for various image sizes.}
\label{fig:runtimes-and-speedup}
\end{figure}

Computing a single  {\pmask} on the CPU takes several seconds no matter 
whether OpenMP or pthreads are used.
When using CUDA and thus the GPU
the total runtimes match the interspike interval constraint. For a typical resolution of $800 \times 600$ for an SLM the computation of~25
iterations in single precision takes 45~ms including transfer of the given {\osp} to the GPU (1~ms).
The iteration essentially converges after one step (Fig.~\ref{fig:results}).
Hence a reasonable {\osp} is available already after much less than 10~ms.
However, the precise figure depends on the problem size and number of iteration steps.
Thus we keep 10~ms as a conservative bound. 
The left panel of Fig.~\ref{fig:runtimes-and-speedup} shows GPU runtimes per iteration in total and broken down into the contributions due to FFT (green and cyan bars) and enforcement of amplitude constraints (blue and grey bars)
for different image sizes and 25 iterations of Eq.~(\ref{e:AP}).
The runtimes are further subdivided into the results for single and double precision indicated by the red and magenta bars. The proportion of work to be done in the FFT increases with problem size as the FFT is of log-linear complexity.
Enforcing the amplitude constraints is linear in the problem size as each pixel is visited only once per iteration and exchange of information between different pixels is not required.

The right panel of Fig.~\ref{fig:runtimes-and-speedup} shows the speedups of the CUDA implementation over its OpenMP and pthreads counterparts. 
On average CUDA is
50 times faster per iteration than the 12-thread CPU variants.
The performance gain per iteration solely depends on the size of an {\osp}. For the Fermi architecture used in the Tesla C2070 the floating point performance in double precision is half of the one for single precision. This is due to the fact that a \verb+double+ is twice as large as a \verb+float+ and thus requires twice as much memory bandwidth.
On CPUs this is less of an issue since they focus on hiding memory latencies by branch prediction.
This is reflected by the fact that for double precision the 
speed\-up is roughly only half of the one for single precision.
Yet, this suffices to get {\osp}s in double precision within the limits set by the interspike intervals.

\mybeginwidetext
~
\myendwidetext

\subsection{Precision and Convergence}

The second issue is the influence of the floating point precision~$\epsilon$ on convergence and performance.
Figure~\ref{fig:results} shows the convergence and reconstruction results for a spot pattern as it would be used in a {\ps} experiment.
Figure~\ref{fig:results-siemens} summarizes the results for the Siemens star which is a standard test image for the spatial resolution achieved by reconstruction algorithms.

The reconstructed {\osp}s are given as inset with a logarithmic gray scale for intensity.
The convergence curves represent the behavior of the physical error,~Eq.~(\ref{eq:phys-error}), 
and the gap,~Eq.~(\ref{eq:gap-error}), with respect to the number of iterated steps in Eq.~(\ref{e:AP}).
We are interested in the influence of the hardware architecture and the precision. Hence the convergence history is given for single and double precision on GPU and CPU.
The details of the curves for the gap and for the physical error in the insets reveal that the behavior primarily depends on whether the computation is run in single or double precision but not on the architecture.
This is a subtle effect on the order of the single precision accuracy 
as illustrated by the scaling of the ordinates in the insets.
Both figures
show that convergence of the {\pmask} in single 
is as good as in double precision as each quantity has a unique limit value independent of the precision. 

The intensity plots of the reconstructed  {\osp}s indicate that the contrast between
lit and unlit areas is 3 orders of magnitude.
A look at the center of the Siemens star shows that structures down to a few pixels can be resolved.
The insets of figures~\ref{fig:results} and~\ref{fig:results-siemens} demonstrate that 
the gap, as defined in~Eq.~(\ref{eq:gap-error}), and the physically motivated error,~Eq.~(\ref{eq:phys-error}), saturate within one iteration indicating the inconsistency of our~{\prp}.
All further changes are~${\cal O}(\epsilon)$.
Convergence does not depend on hardware as the limit values of error and gap are several orders of magnitude larger than any precision. 
Depending on~$\epsilon$ we expect the following resolution limits for~$G$.

Our number of pixels is of the order of~$10^{6}$. 
Assuming statistical independence for the errors of~$u^{\nu}_{jk}$ we get as theoretical limit $G_{th} \propto 10^{3}\epsilon$, i.e. $10^{-5}$ for~single precision~($\epsilon=10^{-8}$) and $10^{-13}$ for~double precision~($\epsilon=10^{-16}$).
An interesting phenomenon reflecting the difference between exact and finite precision arithmetic is revealed by comparing the convergence behavior as function of~$\epsilon$. 
Single precision (blue and red curves) cannot resolve the inconsistency of the {\prp}, i.e. whether or not
$S\cap M=\emptyset$, as~$G\approx G_{th}$. According to~\cite{Luke11} this should improve convergence.
The downside is, that while the {\prp} appears to be consistent from a numerical point of view, larger~$\epsilon$ means worse approximation of the projection operators. For double precision (green and magenta curves) we get more accurate projection operators but at the same time the gap is resolved as for both precisions~$G$ is of similar magnitude.
This renders the {\prp} inconsistent again, justifying our assumption of inconsistency.
Our results also show that, despite a rather large $G$ the {\maproj} does not suffer from stagnation at bad local minima which 
otherwise would call for more sophisticated algorithms like RAAR~\cite{Luke2005}.

\section{Conclusions}
\label{sec:conclusion}

Mathematically, computing a phase-only hologram to create an {\osp} which selects predetermined 
neurons is equivalent to the 
problem of 
wavefront reconstruction.
Useful approximations of a {\pmask} for a given {\osp} can only be achieved by
iterative algorithms like the widely used Method of Alternating Projections.

From the point of view of software engineering we have shown how to integrate CUDA into a complex software environment in a modular way without sacrificing performance. 
The high modularity of our simulation framework has several key advantages.
Code redundancy is minimized.
The template techniques let the code reflect the mathematical structure of the problem.
The effort to switch between the three parallelization techniques tested, CUDA, OpenMP and pthreads, is reduced to a single word in a single line of code and can be done either at compile or at run time. The framework makes it easy to implement other reconstruction algorithms and to apply it to other problems of wavefront reconstruction totally unrelated to the presented test case from the field of optogenetics.
For instance, we could integrate the ideas discussed by Thalhammer et al. for speeding up the switching of liquid crystal SLMs~\cite{Thalhammer:OpEx13}.
Typical switching times are of the order of 10~ms and thus may interfere with the spiking dynamics of the neurons. 

Finally, only the CUDA-based implementation is capable of the necessary frame rates for stimulating networks of
optogenetically altered neurons on their intrinsic timescale of several ms. 
Our results show that at most 5 iterations suffice 
to compute a phase mask within less than 10ms,  matching the time-scale of the dynamics of neural activity. 

\section*{Acknowledgments}
%{\it Acknowledgements}:
SCK thanks NVidia for support under the terms of the CUDA Teaching Center G\"ottingen headed by Gert Lube. 
DRL was supported in part by DFG grant SFB755-C2. 
We thank Hecke Schrobsdorff  and Fred Wolf for pointing us at the problem and fruitful discussions. We acknowledge support by the German Research Foundation and the Open Access Publication Funds of the G\"ottingen University.

\bibliographystyle{abbrv}

\end{document}